\documentclass{article}
\usepackage{jcompbio}
\usepackage{makeidx}
\usepackage{graphicx}
\usepackage{datetime}
\usepackage{amsmath}
\usepackage{appendix}

\begin{document}
\title{Fast and Accurate Detection of Multiple QTL}
\author{Carl Nettelblad*, Behrang Mahjani, Sverker Holmgren\\
\texttt{\{carl.nettelblad, behrang.mahjani, sverker.holmgren\}@it.uu.se} \\
Division of Scientific Computing\\
Department of Information Technology\\
Uppsala University\\
Uppsala, Sweden \\
Phone: +46 18 471 00 00}
\maketitle
\abstract

We present a new computational scheme that enables efficient and reliable Quantitative Trait Loci (QTL) scans for experimental populations. Using a standard brute-force exhaustive search effectively prohibits accurate QTL scans involving more than two loci to be performed in practice, at least if permutation testing is used to determine significance. Some more elaborate global optimization approaches, e.g.\ DIRECT, have earlier been adopted to QTL search problems. Dramatic speedups have been reported for high-dimensional scans. However, since a heuristic termination criterion must be used in these types of algorithms the accuracy of the optimization process cannot be guaranteed. Indeed, earlier results show that a small bias in the significance thresholds is sometimes introduced.

Our new optimization scheme, PruneDIRECT, is based on an analysis leading to a computable (Lipschitz) bound on the slope of a transformed objective function. The bound is derived for both infinite size and finite size populations. Introducing a Lipschitz bound in DIRECT leads to an algorithm related to classical Lipschitz optimization. Regions in the search space can be permanently excluded (pruned) during the optimization process. Heuristic termination criteria can thus be avoided. Hence, PruneDIRECT has a well-defined error bound and can in practice be guaranteed to be equivalent to a corresponding exhaustive search. We present simulation results that show that for simultaneous mapping of three QTL using permutation testing, PruneDIRECT is typically more than 50 times faster than exhaustive search. The speedup is higher for stronger QTL. This could be used to quickly detect strong candidate eQTL networks.

\newpage
\section{Introduction}
The rapid development of off-the-shelf technology for molecular genetics has resulted in that dense genetic maps and the corresponding genotype information can be provided in a much easier and cheaper way than before. This development opens new possibilities for analysis of quantitative traits, i.e. traits that exhibit a continuous phenotype distribution. Since most important traits in humans, animals and plants can indeed be seen as quantitative and affected both by the genetic composition and the environment, genetic mapping of them represents both a major opportunity and a challenge for modern genetics.

The underlying genetic architecture of a quantitative trait can be described by identifying a set of \emph{quantitative trait loci} (QTL) in the genome for a population and attributing effect values to these loci using a suitable statistical model framework. The standard approach for locating a QTL is based on \emph{interval mapping} (IM) \citep{IM}. Evaluating the standard IM model at a given position in the genome involves solving a maximum-likelihood problem based on genotype and phenotype frequencies for the population studied. In a \emph{QTL search}, the evaluation of the statistical model is repeated for a large set of candidate positions in the genome to determine the QTL locations that results in the best model fit. Mathematically, this corresponds to solving a global optimization problem using some optimization scheme.

The result of a QTL mapping procedure is normally useful only if a proper significance threshold can be derived. Already when searching for a single QTL, traditional $\chi^2$ approximations have been shown to have a significant bias \citep{carbonell}. Therefore, randomization testing is frequently used \citep{Churchill1994}, where many permuted datasets and the corresponding QTL searches are employed to empirically derive the distribution of the optimal model fit under the null hypothesis of no QTL being present.

In general, it can be assumed that multiple interacting QTL (epistatic interactions) should be included in a model to fully describe the genetic effect on a trait \citep{doerge2002}. However, using a $d$-QTL model with general interactions results in that a $d$-dimensional global optimization problem has to be solved for each QTL search. For QTL models the optimization landscape is often varied with many local minima, and still today the standard approach for QTL mapping problems is to use a brute-force exhaustive search over a dense lattice covering the search space. For multidimensional searches this approach rapidly becomes computationally intractable. Beyond $d=2$, QTL mapping employing true multidimensional dimensional optimization has, as far as we know, not been used in practice. This means that geneticists have so far not drawn any firm conclusions on how important epistatic interactions between more than two QTL are for describing quantitative traits. However, there are indications that such interactions may indeed be important, see e.g.\ \citep{Haley04}.

A main reason for the inconclusive situation regarding the importance of more complex epistatic interactions is the lack of efficient and reliable computational tools for performing multi-dimensional QTL mapping, as well as determining the joint significance of the set of QTL found. In this context, it should also be noted that even under the assumption that the QTL do not interact (i.e. the true population effects are perfectly additive), the estimated QTL effects will be more accurately estimated if all putative loci are included in a single, multi-locus model. Hence an efficient computational tool can be very useful also for such settings. If a single-locus model is used repeatedly, correlations between loci could inflate or distort the detected effects. This is especially true for linked loci residing on the same chromosome, but also holds due to random corelations for loci located on different chromosomes.

Some early examples of simultaneous mapping of two QTL are found in \citep{Carlborg2000}, where a genetic optimization algorithm is used for solving the optimization problem. Today, performing two-dimensional QTL mapping is regarded as a standard procedure by many researchers in genetics. In \citep{kl3} the deterministic optimization algorithm DIRECT \citep{Jones1993} was introduced for solving QTL search problems also for $d\geq 2$, but these results have so far not been used to perform high-dimensional QTL mapping experiments of relevance to genetics.

In this paper, we present a computational scheme that enables efficient and reliable solution of QTL mapping problems in experimental populations for high dimensional general models of interacting QTL. Our scheme is based on an analysis of the behavior of the objective function (the results from a linear regression QTL model fit) and implemented in deterministic global optimization framework. We also show how the result of the analysis and our optimization framework can be used to set up permutation testing in a very efficient way to determine the relevant significance thresholds. Our algorithms are structured in such a way that the problem of determining the set of $d$ QTL resulting in the optimal model fit is separated from the evaluation of the model of the QTL effects. This implies that models based on e.g.\ both the linear regression approximation and the standard interval mapping maximum-likelihood model can easily be included in a production tool for genetic analysis. We focus on linear regression models in this paper since they are much less computationally demanding and lend themselves to the type of transformations that are exploited in our analysis of the behavior of the objective function. In this context it is important to note that if complete genetic information is available, and general assumptions of normal distributions hold, the linear regression and maximum-likelihood are equivalent \citep{hk92}.

\section{Linear Regression QTL Models}
We consider QTL analysis for experimental populations with known relations between individuals, with the founder individuals demonstrating some origin-defining feature. This can be a matter of a known genetic relation (a set of common inbred or outbred lines) between the founders, or founders expressing a specific phenotype. Today, very dense marker maps are available and we formulate the analysis assuming information on allele origin being available in any position tested.

If a model with a total of $d$ QTL is used, then the corresponding $d$ sites in the genome are assumed to represent the only genetic factors that contribute to the phenotype. One can split individuals based on genotype into $2^d$ classes in a backcross ($3^d$ for an intercross), since each of the $d$ sites can take $2$ different values ($3$ for an intercross, two homozygote genotypes and one heterozygote, if allele parental origin is ignored). Within each of these classes the variance is entirely non-genetic. Another important assumption is that phenotypes are normally distributed with different means, but identical variance in each class, i.e. $y|g_1,g_2,...g_d \sim N(\mu_{g_1,g_2,...g_d},\sigma^2)$. For details on different experimental cross structures see e.g.\ \citep{statisticsb}. For ease of presentation, we now consider the typical case where we have an $F_0$ generation of individuals that can be considered to belong to either out of two lines, $Q$ and $q$. Assuming loci act additively, one can model a relation between genotype and phenotype in a backcross based on these founders as:
\begin{align}
y_i=\mu(x)+\sum_{j=1}^d a_j(x_j)Z_{i}(x_j)+e_i(x)
\text{ , where: }Z_{i}(x_j)=\left\{
       \begin{array}{ll}
         1, & \hbox{Qq at $x_j$;} \\
         0, & \hbox{qq at $x_j$.}
       \end{array}
     \right.\label{reg1}
\end{align}
In this model $x=(x_1,...,x_d)$ is a vector of $d$ elements which defines the search space, spanning a hypercube where $x_i$ ranges over the length of the genome. In practice, the search space volume can be slightly reduced by employing the symmetries resulting from that the ordering of the QTL within the model is irrelevant. The phenotypes for all individuals are denoted by $y=(y_1,...,y_n)$ denotes a vector of $n$ elements, $e_i$ is normally distributed with mean $0$ and variance $\sigma^2$, $\mu$ is the reference effect and $a$ is the additive effect. The model (\ref{reg1}) can be written in a matrix form:
\begin{align}
y=\textbf{A}(x)b(x)+e(x)\text{, where } b(x)=(\mu(x),a(x))^T. \label{reg1matrix}
\end{align}
The least-squares estimate of the QTL effects for this linear model is:
\begin{align}
\hat{b}=(A(x)^{T}A(x))^{-1}(A(x)^Ty)\\
\hat{\sigma}^2=\frac{1}{n}(y-A(x)b(x))^T(y-A(x)b(x)).
\end{align}
These QTL positions can now be found by minimizing the residual sum of squares over all $x$ and $b$:
\begin{align}
RSS_{opt}=\underset{x,b}{\operatorname{min}}(y-A(x)b(x))^T(y-A(x)b(x)).\label{fullopt}
\end{align}
The solution to this minimization problem can be separated into two parts; The inner, linear problem:
\begin{align}
RSS(x)=\underset{b}{\operatorname{min}}(y-A(x)b(x))^T(y-A(x)b(x)),
\end{align}
and the outer, non-linear problem:
\begin{align}
RSS(x)=\underset{x}{\operatorname{min}}RSS(x).
\end{align}
Solving the minimization problems (\ref{fullopt}) for a multiple QTL mapping problem is computationally heavy since $x$ is a d-dimensional vector and the optimization landscape for the outer (global) optimization problem is in general quite complex. It is clear that an optimal (albeit not necessarily unique) solution to the QTL search problem always exists, but to determine if a result is statistically valid a significance threshold has to be determined. If permutation testing is used for determining the genome-wide threshold, several hundreds or thousands of QTL mapping problems of the type (\ref{fullopt}) must then be solved for the permuted datasets.

\section{A global optimization algorithm for QTL mapping}\label{theory}
Our new computational scheme for solving the high-dimensional QTL mapping problem is based on a deterministic Lipschitz optimization approach implemented in the DIRECT \citep{Jones1993} algorithmic framework. In this section we first review the original DIRECT algorithm and the summarize how this scheme was adapted to solve multidimensional QTL search problems in \citep{kl3}. We then present the basic idea behind our new optimization scheme, named \emph{PruneDIRECT} as a background for the more detailed description and analysis in later sections.

\subsection{The original DIRECT algorithm}
In the DIRECT scheme, the search space is successively divided into progressively smaller boxes. The search effort is focused in the most promising regions and the subdivision of less promising boxes is postponed. For a traditional exhaustive search over a $d$-dimensional hypercube, the objective function is evaluated in a brute-force fashion at every point of a fine lattice covering the search space, using no further information on the objective function. If DIRECT is run to completion using a minimum resolution criterion matching the step length in the exhaustive search lattice, it will eventually explore exactly the same points as the corresponding exhaustive search. The efficiency of the DIRECT algorithm comes from exploiting the property that the most promising regions are explored first. A heuristic termination criterion is then needed to stop the search well before the it devolves into an exhaustive search. As the heuristic criterion has no firm mathematical foundation, this implies that there is no well-established guarantee that the result from DIRECT is equal to the result from an exhaustive search over the same space. This is a typical result for global optimization schemes aimed at solving general problems. To be able to provide guaranteed accuracy without exploring the full search space, more information on the properties of the objective functions must be provided and used in the optimization procedure.

\begin{figure}
\centerline{
\includegraphics[trim = 0mm 0mm 0mm 0mm, width=7cm]{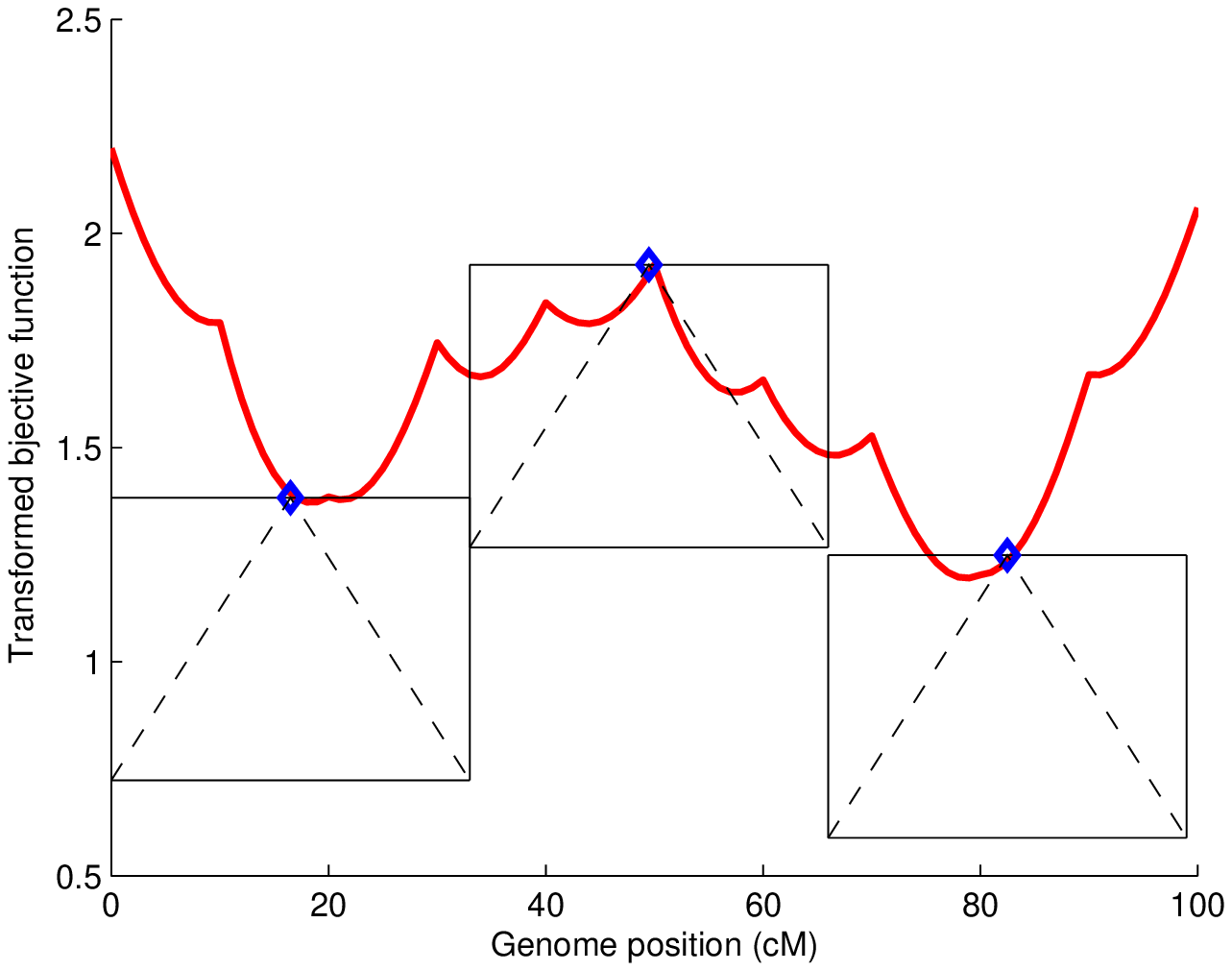}
\includegraphics[trim = 0mm 0mm 0mm 0mm, width=7cm]{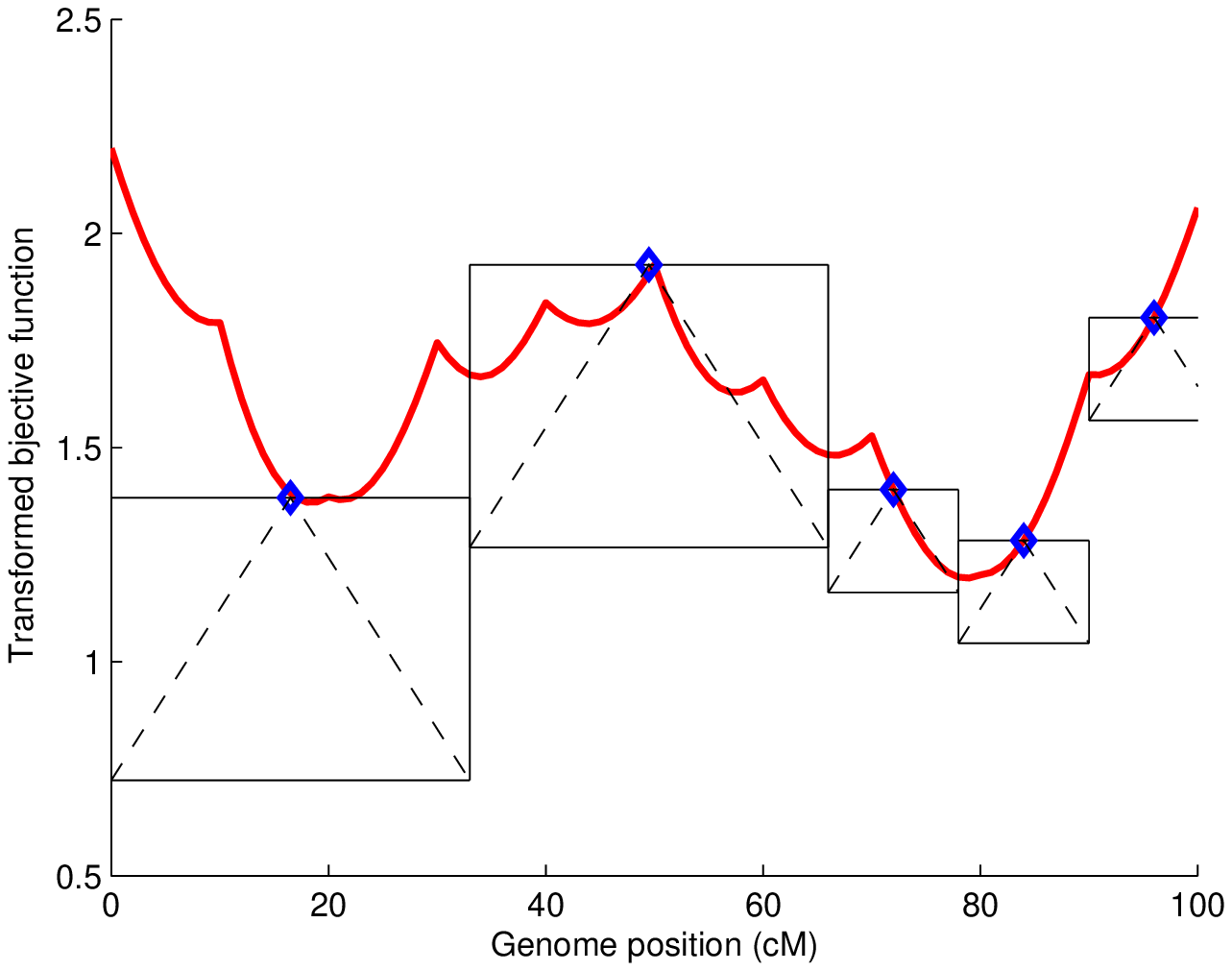}
}
\caption{\footnotesize\emph{To the left, three boxes, in the optimization space. As all boxes are of equal size, only one will be selected for splitting in the convex hull, resulting in the boxes to the right. If the splitting would continue, two boxes would be split, one from each size, as the smallest function value in the smaller box size is slightly lower than the smallest value in the larger one. Dashed lines indicate possible minimum function values found at each distance from the box centroid, assuming a strict Lipschitz bound of $K = 0.04$.}}\label{shadedarea}
\end{figure}

The original DIRECT algorithm initially creates a single search box covering the full search space. The objective function, i.e. in our case the RSS, is evaluated in the center of this box. The box is then split into three equally sized boxes along the majoring dimension. This trinary split results in the centroid of the original box coinciding with the centroid of one of the new boxes. Therefore, only two additional function evaluations are required for the three resulting boxes. DIRECT then continues iteratively splitting the boxes. At the end of each DIRECT iteration, the convex hull is determined among the remaining boxes, in a space of box radii vs. objective function values. This hull determines which boxes to split in the following iteration. By computing the hull, the RSS value for the sequence of boxes will be monotonously decreasing if tracing along the box radii. The idea is that a promising box is characterized by either being large, so there is a high possibility that exploring it further might recover a new optimum, or that the value at the centroid is close to the current optimum and hence a new optimum might be found, even if the radius is small. This qualitative argument can also be represented as different assumptions of the value of a Lipschitz constant $K$, i.e. on the maximal slope, of the objective function \citep{Jones1993}. Figure \ref{shadedarea} illustrates a few iterations of DIRECT in a simple one-dimensional space. The hull is ``peeled off''
when those boxes are split, making new boxes available in the next
iteration. This process is repeated until a suitable termination
condition has been reached.

\subsection{Adaption of DIRECT to solve QTL mapping problems}
In \citep{kl3}, the original DIRECT algorithm was adapted for solving QTL search problems. A heuristic termination criterion, as described above, was used to stop the search and the resulting scheme was proven to be orders of magnitudes faster than the corresponding exhaustive search. This made QTL searches for at least $d=3,4$ and $5$ possible using a reasonable computational effort. However, the results in \citep{kl3} also show that using DIRECT for the permutation tests can result in a bias in the significance thresholds compared to those from an exhaustive search. The reason for this is that the termination criterion used in \citep{kl3} results in premature termination for the very flat optimization landscape present in most permuted cases. The end result is that a putative $95 \%$ threshold instead would give e.g.\ $94.8 \%$ significance. It should also be noted that no bias was detected in the DIRECT searches for non-permuted phenotype data, probably due to the more structured nature of these optimization landscapes. A main result of the analysis presented in later sections is that it is possible to terminate our new PruneDIRECT process at a much earlier point than the corresponding exhaustive search while still being able to \emph{guarantee}, up to some defined residual probability, that the same global optimum is found as for the exhaustive search.

We argue that the successful results for DIRECT presented in \citep{kl3} originate from the fact that the probability of co-inheritance of two genomic loci on a single chromosome is related to their physical distance from each other. The standard unit for genetic distances, mapping distance, is even directly defined from such probabilities. In \citep{kl3} this fact is not explicitly used in the algorithm, but since the end result is that the slope of the objective function is limited, the performance of DIRECT ends up as being quite good. It is also clear that this notion of co-inheritance cannot be extended across chromosome boundaries. Separate chromosomes segregate independently and can be considered to have an infinite mapping distance. In the search space for the QTL mapping problems all genome locations are lined up as a single line, one chromosome after another. There is a clear disjunction at the chromosome borders, corresponding to a discontinuity of the RSS at these boundaries. In \citep{kl3}, this situation is handled by introducing what is essentially several DIRECT searches governed by a common priority queue. In this version of the algorithm, each chromosome combination is considered to correspond to a separate search space, called a chromosome combination box (cc-box). In the first step of the DIRECT scheme, the RSS is evaluated at the centroids of all these cc-boxes.

\subsection{The PruneDIRECT algorithm}
\label{prunedirect}
The original DIRECT algorithm is based on an assumption of Lipschitz continuity in the objective function, but with an unknown Lipschitz constant. This means that we know that there exists a constant $K$ such that no partial first-order derivative of the RSS will exceed $K$ in any position, but that the value on this bound of the slope of the objective function is not known.
Our improvement of the DIRECT procedure for QTL mapping is based on that a Lipschitz bound for a QTL problem can be computed, based on the relation between co-inheritance and physical distance mentioned above. Parts of the search space can then be permanently excluded from further subdivisions at each iteration in the DIRECT procedure, resulting in an efficient algorithm where a heuristic termination criterion does not have to be used. The resulting scheme has a well-defined error bound and is equivalent to the corresponding exhaustive search.

If a bound on the Lipschitz constant $K$ is known it is possible to compute upper and lower bounds for the objective function within any box in the search space, given knowledge of the value at the centroid. When using the original DIRECT algorithm, $K$ is unknown but it is still possible to impose a partial ordering of boxes; if a box $A$ is both smaller and has a larger value of RSS than another box $B$, then no value of $K$ can result in a smaller minimum bound within $A$ than within $B$. This partial ordering gives rise to the choice of a convex hull of boxes to split in the next iteration.

However, DIRECT can also be considered in the more general context of Lipschitz optimization schemes. Here, the more traditional schemes \citep{lo0,lo1,lo2,lo3} assume that a bound of the Lipschitz constant is known. By introducing an upper bound on $K$ in DIRECT, the resulting algorithm is brought close to classic Lipschitz optimization. A known value of the objective function in small boxes will also enable exclusion of larger boxes, hereby introducing a pruning of the search tree, indicating the choice of name for our new algorithm. In contrast to the peeling effect of successive convex hulls in DIRECT, which will eventually result in every existing box getting split, box volumes pruned due to a Lipschitz criterion are permanently irrelevant for further evaluation of the RSS and can be removed from all data structures.

The natural termination condition for the PruneDIRECT algorithm is given by the finite resolution criterion corresponding to the lattice in the underlying exhaustive search. This results in an error bound of $Kh$ on the value of the objective function, where $h$ is the step length in the lattice. Thus, a bound on the Lipschitz constant $K$ leads to that the new PruneDIRECT algorithm has improved performance compared to exhaustive search (by excluding parts of the search space) and a well-defined error bound. Since no heuristic termination criterion is used anymore, a main new result is also that PruneDIRECT can be guaranteed to be equivalent to the corresponding exhaustive search, effectively removing any bias of the results caused by the optimization procedure.

For multi-dimensional QTL mapping that explore the full search space it is possible to directly use the permutation test methodology for single-QTL models. However, performing hundreds or thousands of multi-dimensional QTL searches for permuted data in order to get a significance threshold can be a very computationally demanding task, even when an efficient global optimization scheme is employed. Here, it should be noted that since the purpose of doing the random permutations is to determine a significance threshold empirically, it is not necessary to locate the location of the best model fit for every set of randomized data. Rather, it is enough to answer the yes/no-question if it possible to find a location in the permuted dataset with a residual variance below that of the putative set of QTL. If the optimum value used to determine pruning is replaced by the value of the QTL candidate tested, rather than the optimum found in the current permute search, significant decreases in computational effort are possible, even compared to using the same algorithm to solve the full QTL search problem for each permuted dataset. The significance levels derived are identical. We use this approach to derive a very fast scheme for performing permutation tests, described in more detail later.

\section{Presenting the LogVar objective function and its Lipschitz bound}
In this section, the behavior of a transformed QTL search objective function at, and in the vicinity of, a QTL is considered. For a QTL search, the explainable genetic variance can be considered a natural objective function since the position with minimum residual variance can be defined as the location of a putative QTL. The benefit of using our transformed objective function, which we call LogVar, is the possibility to calculate a bound for its Lipschitz constant. The derivation of the bound is based on calculating the explainable variance as a function of genetic distance from the true QTL. We start by presenting the transformation and deriving the bound for an infinite size population and then move to the more realistic situation where QTL mapping for a finite size population is performed.

\subsection{Infinite size population}
\label{expvarsect}
Consider an infinite size population. The total variance is the sum of genetic variance and environmental variance. If there are $d$ QTL, then all genetic variance is explainable by them. We start by analyzing $d=1$ and then show how these results can be generalized to higher dimensions.

Assume there is a QTL at position $x_0$. The total variance is the "mean squared error". We split the total variance into a sum of undiscovered genetic variance and the discovered genetic variance at position $x+x_0$:
\begin{align}
\nonumber V_{total}=V_{genetic}+V_{enviromental}=V_{g-discovered}(x+x_0)+V_{g-undiscovered}(x+x_0)+V_{e},\\
\text{which can be abbreviated into } V_{t}=V_{gd}(x+x_0)+V_{gu}(x+x_0)+V_{e}.
\end{align}

Then define $V_r(x+x_0)=V_{gu}(x+x_0)+V_e$ as the unexplainable variance at position $x+x_0$. The goal is to express $V_r(x+x_0)$ as a function of the recombination frequency and then find a bound for it. For simplicity, we start by letting $V_e=0$ and later we add $V_e$ back to the calculation.

Assume that the phenotype values for the two QTL genotypes are $0$ and $1$. Due to the symmetric structure of individuals with phenotype value $0$ and individuals with phenotype value $1$, $V_r(x+x_0)$ for each class is equal to the total residual variance. Hence:
\begin{align}
\nonumber V_r(x+x_0)=V_r(x+x_0 \text{ for class 1})&=\text{Prob(no recombination)}&\times (1-\mu(x+x_0\text{ for class 1}))^2\\
&+\text{Prob(recombination)}&\times (0-\mu(x+x_0\text{ for class 1}))^2 \label{vr1}
\end{align}
where:
\begin{align}
\nonumber \mu(x+x_0 \text{ for class 1})&=1\times \text{Prob(No recombination)}+0\times \text{Prob(recombination)}\\
 &= p(x+x_0).\label{mu}
\end{align}
Hence, from substituting (\ref{mu}) in (\ref{vr1}) and denoting the recombination frequency at position $x+x_0$ by $p(x+x_0)$ we get:
\begin{align}
\nonumber V_r(x+x_0)&=(1-p(x+x_0))(1-p(x+x_0))^2+p(x+x_0)(0-p(x+x_0))^2\\
&=p(x+x_0) -p(x+x_0)^2. \label{eq.7}
\end{align}
This result for $V_r$ should be related to the total phenotypic variance, which is the genetic variance at position $x_0$, since we defined $V_e = 0$. We know that all the variance at a QTL point is the discovered variance, hence:
\begin{align}
V_t=V_g(x_0)=V_{gd}(x_0)+v_{gu}(x_0)=[0.5\times (1-0.5)^2+0.5\times (0-0.5)^2]+0=0.25. \label{eq.8}
\end{align}
At last, we need to have a recombination map, relating genetic distance to the recombination frequency. We use Haldane's mapping function \citep{haldane}:
\begin{align}
p(x+x_0)=0.5+0.5e^{-\frac{2x}{100}},\label{eq.9}
\end{align}
where $x$ is the genetic distance from a fixed point $x_0$ measured in centimorgan.\\
Inserting $p(x+x_0)$ from (\ref{eq.9}) into equation (\ref{eq.7}) we get:
\begin{align}
V_r(x+x_0)&=(0.5+0.5e^{-\frac{2x}{100}})(1-(0.5+0.5e^{-\frac{2x}{100}}))\\
&=0.25-0.25e^{-\frac{4x}{100}}\\
&= V_t-V_g(x_0)e^{-\frac{4x}{100}}(\text{from \ref{eq.8}}). \label{vrdef}
\end{align}
Here, $V_e$ can be added to the final calculation since variance is additive. Alternatively, it is safe to assume that $V_e$ is already included in calculating $V_t$. Now that we have the formula for calculating the explainable variance as a function of genetic distance, we can introduce a transformation of the objective function in the DIRECT optimization. Ignoring a scaling factor, the RSS is equivalent to the residual variance $V_r$. Adding a constant will not affect the location of minima, so we can instead consider $g(x) = V_r(x) - V_t = -V_g e^{-\frac{|4x|}{100}}$ based on (\ref{vrdef}), assuming $x_0=0$. Furthermore, the function $g(x)$ is always negative, so $f(x) = -\ln (-g(x))$ will always be defined and share the locations of minima with $g(x)$. Hence, we can introduce a transformed objective function, call it LogVar objective function:
\begin{align}
f(x)=-\ln(-V_r(x) + V_t)=-\ln(V_g e^{-\frac{|4x|}{100}}). \label{infinite}
\end{align}
Here, it is easy to verify the derivative of the transformed objective function is bounded:
\begin{equation}
\label{boundedd}
\frac{df(x)}{dx} = \frac{d}{dx} -\ln (V_g e^{-\frac{|4x|}{100}}) = \pm 0.04.
\end{equation}

The use of $|x|$ in the definition of $g(x)$ is related to the fact that
$x$ is defined as the distance from the QTL, while a position in the chromosome can naturally be both upstream and downstream from this position. It is possible to shift the function by introducing the true QTL position $y$, resulting in $f(x) = \ln V_g - 0.04 |x - y|$.

The result above can be further generalized, maintaining the bound on $K$. First, we have assumed that all genetic variance was attributable to a single locus (at $y$). We can now assume a single-locus model for
analysis, but that the true QTL, with respective components of $V_g$
are represented as a vector $\mathbf{y} = \left \{y_1, y_2, ..., y_n
\right \}$, resulting in the following expression for the residual
variance (assuming all $y_i$ being unlinked):

\begin{equation}
\label{multiunlinked}
V_r(x) = V_t - \sum_{i=1}^n V_{g_i} e^{-\frac{4|x - y_i|}{100}}.
\end{equation}

If all QTL are indeed unlinked, the positions $y_i$ relative to any
point considered in a single-QTL model will be $\pm \infty$, except for at most one $y_j$. Since linkage is transitive, the observation position $x$ can at most be linked to a single QTL. Thus,
(\ref{multiunlinked}) reduces to:

\begin{equation}
\label{redmultiunlinked}
V_r(x) = V_t - V_{g_j} e^{-\frac{4|x - y_j|}{100}},
\end{equation}
and the derivative bound on (\ref{redmultiunlinked}) follows from the
result in (\ref{boundedd}).

The next extension is to make the search landscape itself multidimensional, replacing the scalar $x$ with a vector $\mathbf{x} =
\left\{ x_1, x_2..., x_m \right \}$. Assume that each $x_i$ is defined
with a point of reference in linkage with the corresponding QTL
$z_i$. Modeled locations not in linkage with any true QTL will result
in no explainable variance, and therefore do not need to be
considered. As the numbering of both vectors is essentially arbitrary,
all other cases are also symmetrical to this one. Using the earlier
result, each $x_i$ will only have a term for the corresponding $z_i$,
as all other mapping distances $|x_i - z_j|$ would be $+\infty$.
This results in:
\begin{equation}
\label{multimulti}
V_r(\mathbf{x}) = V_t- \sum_{i=1}^n V_{g_i} e^{-\frac{4|x_i - z_i|}{100}}.
\end{equation}
We then have:
\begin{equation}
\label{fmulti}
f(\mathbf{x}) = - \ln (-V_r(\mathbf{x}) + V_t),
\end{equation}

\noindent where the logarithm does not affect the location of minima.

Equation (\ref{multimulti}) can be rewritten as:
\begin{equation}
\label{VrC}
V_r(\mathbf{x}) = V_t - V_{g_j} e^{-\frac{4|x_j - y_j|}{100}} -
\underbrace{\sum_{i=1, i \ne j}^n V_{g_i} e^{-\frac{4|x_i -
      y_i|}{100}}}_{C},
\end{equation}
\noindent for any arbitrary $j \in {N}, j \leq n$. Based on
(\ref{VrC}), (\ref{fmulti}) becomes:
\begin{equation}
f(\mathbf{x}) = - \ln (V_{g_j} e^{-\frac{4|x_j - y_j|}{100}} + C).
\end{equation}

From basic calculus, we know that $|\frac{d}{dx} \ln (e^{kx} + C)| =
|\frac{k e^{kx}}{e^{kx} + C}| \leq \frac{d}{dx} \ln (e^{kx}) = \frac{k
  e^kx}{e^kx} = k$ for any $C \in {R}, C \geq 0$. Hence, all
partial derivatives
$\frac{\partial f  }{\partial x_j}$ are confined by the bound given for
the single-dimensional first derivative presented earlier. Thus, the LogVar objective function
$f(\mathbf{x})$ as defined above will again have a well-defined Lipschitz
bound.

\subsection{Explainable Variance as a Function of Genetic Distance (finite size population)}

The bound derived for an infinite size population does not apply directly to the derivative of the
actual residual variance in experimental data, but to the expected value of the
derivative, corresponding to the relation between the mapping distance
and the expected number of crossover events. Depending on what
recombinations are actually present (i.e. in which individuals, with
accompanying phenotype values), the actual residual variance can, and
will, be different from the relationship predicted. This is an effect of
sampling, which should decrease with an increasing size in population
and vanish at a theoretical infinite population size. However, experimental
populations tend to be rather small, and thus the infinite-size approximation cannot
be used directly.

In this section, we derive an approximation to the distribution of the residual sum of errors in the vicinity of any point in the search space for a backcross. We revisit the linear model for a single QTL. Assuming there is a QTL at point $x_0=0$, the model is presented as:
\begin{align}
y_i=\mu(x)+a(x)Z_{i}(x)+e_i(x)
\text{ , where: }Z_{i}(x)=\left\{
       \begin{array}{ll}
         1, & \hbox{Qq at $x$;} \\
         0, & \hbox{qq at $x$.}
       \end{array}
     \right.\label{reg1}
\end{align}
Our goal is to calculate the distribution of residual sum of squares (RSS) at point $x$. At first, we introduce some definitions:
\begin{align}
n&=\text{ total number of individuals}\\
n_0(x)&=\text{ number of 0's at $x$}\\
n_1(x)&=\text{ number of 1's at $x$}\\
m_{10}(x)&=\text{ number of 0's at $x$ which are $1$ at zero}\\
m_{01}(x)&=\text{ number of 1's at $x$ which were $0$ at zero}\\
M_{10}(x)&=\text{ set of indices of individuals which recombined from $1$ at zero to $0$ at $x$}\\
M_{01}(x)&=\text{ set of indices of individuals which recombined from $0$ at zero to $1$ at $x$}\\
p&=\text{ recombination frequency}\\
S_0 &= S_1 = \text{within-class environmental variance} = 1
\end{align}
The last definitions are simplifying assumptions relating to the normal residual assumption underlying linear regression. Based on the above definitions we get:
\begin{align}
n_0(x)&=n_0(0)+m_{10}(x)-m_{01}(x)\\
n_1(x)&=n_1(0)+m_{01}(x)-m_{10}(x)\\
\overline{Z}(x)&=\frac{n_1(x)}{n}
\end{align}
For now, consider $m_{01}(x)$ and $m_{10}(x)$ to be fixed and known. We also assume that $n_0(0)=n_1(0)=n/2$ for simplicity. We introduce the estimates for $\hat\mu(x)$, $\hat a(x)$ and $RSS_x$. Set $a_1=\sum (Z_i(x)-\overline{Z}(x))y_i$ (proportional to the covariance between $Z$ and $y$) and $a_2=\sum(Z_i(x)-\overline{Z}(x))^2$ (proportional to the variance of $Z$), then:
\begin{align}
\hat a(x)&=\frac{a_1}{a_2} \label{ax}\\
\hat\mu(x)&=\overline{y}-\hat{a}(x)\overline{Z}(x) \label{mux}
\end{align}
After some calculation we get the values for $a_1$ and $a_2$:
\begin{align}
a_1&= \sum_{i\in N\setminus N_r(x)}(Z_i(x)-\overline{Z}(x))y_i + \sum_{i\in N_r(x)}(Z_i(x)-\overline{Z}(x))y_i\\
&=\sum_{i\in N_r(x)}(1-2Z_i(0))y_i+\sum_{i\in N} (Z_i(0)-\overline{Z}(x))y_i \label{a1}\\
\nonumber \text{ where }& N=\{1,2,...,n\}, \text{ and }N_r(x)=\{\text{indices for recombined individuals at point x}\}
\end{align}
As equation (\ref{a1}) shows, we split $a_1$ into two sums, one of recombinants and the other non-recombinant, which we can be summarized as:
\begin{align}
a_1&=a_{11}+a_{12}\\
\nonumber  \text{where:}&\\
  a_{11}&=\sum_{i\in N_r(x)}(1-2Z_i(0))y_i=\sum_{i\in M_{01}}y_i-\sum_{i\in M_{10}}y_i \label{a11}\\
  a_{12}&=\;\;\:\sum_{i\in N} (Z_i(0)-\overline{Z}(x))y_i
\end{align}
 $a_{11}$ is a random weighted sum of phonotypes of recombinants, and is the only stochastic variable involved in $\hat a(x)$. The stochastic behavior of $a_{11}$ comes from the fact that $N_r(x)$ is a random variable related to the recombination process. Generally speaking, one can say that $a_{11}$ captures all the stochasticity of $\hat a(x)$. We approximate the value of $a_{11}$ by assuming $y_i$'s in $a_{11}$ to be normally distributed:
 \begin{align}
\text{Assume (only for recombinant individuals): }&y_{i}\sim N(\mu(0)+a(0)Z_i(0),1)
\end{align}
then:
\begin{align}
 a_{11}&\sim N(\mu_{a_{11}},\sigma^2_{a_{11}})\\
\nonumber \text{where from (\ref{a11}):} &\\
\mu_{a_{11}}&=\mu(0)(m_{01}(x)-m_{10}(x))-a(0)m_{10}(x),\\
\sigma^2_{a_{11}}&=S_0m_{01}(x)\bigg(1-\frac{m_{01}(x)-1}{n_0(x)-1}\bigg)+S_1m_{10}(x)\bigg(1-\frac{m_{10}(x)-1}{n_1(x)-1}\bigg)\\
S_0&=1,S_1=1
\end{align}
We are using population correction coefficients $(1-\frac{m_{01}(x)-1}{n_0(x)-1})$ and $(1-\frac{m_{10}(x)-1}{n_1(x)-1})$, since we are sampling from a finite size population. If e.g.\ all individuals would recombine, all randomness in $a_{11}$ would also disappear, which this correction reflects.

We also need to simplify $a_{2}$. One gets after some calculations:
\begin{align}
a_2&=n_0(x)\overline{Z}(x)^2+n_1(x)(1-\overline{Z}(x))^2
\end{align}

Now that we have $a_1$ and $a_2$, we should calculate $RSS$. If a dataset is divided into disjoint categories, it is known that one can write $RSS_x$ as:
\begin{align}
RSS_x=\sum(y_i-\overline{y})^2-n_0(x)(\hat\mu_0(x)-\overline{y})^2-n_1(x)(\hat\mu_1(x)-\overline{y})^2
\end{align}
where $\hat\mu_0(x)$ and $\hat\mu_1(x)$ are the means of each category. Using this formula we get:
\begin{align}
\nonumber RSS_x&=RSS_{total}-n_0(x)(\hat\mu_0(x)-\overline{y})^2-n_1(x)(\hat\mu_1(x)+\hat a(x)-\overline{y})^2 \\
&=RSS_t-a_2\hat a(x)^2 \label{RSSx}
\end{align}

Given the above description, we calculate the cumulative distribution function (CDF) of $RSS_x$:
\begin{align}
\nonumber F_{RSS_x}(rss_x)&=P_{a_{11}} \left (\hat a(x)^2 \leq \frac{rss_x-RSS_t}{-a_{2}} \right )\\
\nonumber &=1&-P_{a_{11}}(a_{11}\leq \sqrt{(RSS_t-rss_x)a_2}-a_{12})\\
& &+P_{a_{11}}(a_{11}\leq -\sqrt{(RSS_t-rss_x)a_2}-a_{12})
\end{align}
Set $R=\sqrt{(RSS_t-rss_x)a_2}$ and then normalize the CDF to get the standard normal:
\begin{align}
\text{Set: } z&=\frac{a_{11}-\mu_{a_{11}}}{\sigma_{a_{11}}}, \text{ then:}
\end{align}
\begin{align}
F_{RSS_x}(rss_x|m_{01}(x),m_{10}(x),y,Z(0),p)=1-\Phi \left (\frac{R-a_{12}-\mu_{a_{11}}}{\sigma_{a_{11}}} \right )+
\Phi \left (\frac{-R-a_{12}-\mu_{a_{11}}}{\sigma_{a_{11}}} \right )\label{Phi}
\end{align}
where $\Phi$ is the CDF of standard normal distribution. Given this result, we remove the constrains on $m_{01}(x)$ and $m_{10}(x)$ by summing over all values they can take. Since $n_0(0)=n_1(0)=n/2$, we have $m_{01}(x),m_{01}(x)\in \{0,1,...,n/2\}$. Hence, one gets:
\begin{align}
&F_{RSS_x}(rss_x|y,Z(0),p)= \nonumber\\
&\label{frssx}\sum_{m_{01}(x)} [B(m_{01}(x),n_0(0),p) \sum_{m_{10}(x)} B(m_{10}(x),n_1(0),p)F_{RSS_x}(rss_x|m_{01}(x),m_{10}(x),y,Z(0),p)]
\end{align}
Where $B$ is the binomial probability mass function, and the sums are over all values of $m_{01}(x),m_{01}(x)\in \{0,1,...,n_{0/1}(0)\}$.

In this section, we derived the CDF of $RSS_x(x+x_0)$, the residual sum of squares at point $x$ in the vicinity of $x_0$, conditional on known $m_{01}$ and $m_{10}$, by approximating the phenotypes of recombinants with a normal distribution. Then, we removed the condition on $m_{01}$  and $m_{10}$ by computing a probability-weighted sum over all values they can take. This sum forms an approximation of the CDF for $RSS_x$ for fixed $y_i$'s. This CDF is then a two-levelled binomially weighted sum of mixture normals. One can find the value of $rss_x$ for which this CDF is arbitrarily close to $1$, i.e. the $1-\epsilon$ quantile. We thus have an upper bound for the value of $rss_x$ at location $x + x_0$, which holds with probability $1-\epsilon$.

\subsection{Applications} \label{appl}
Section \ref{prunedirect} for an infinite-size population, the PruneDIRECT algorithm and
the Lipschitz bound for the LogVar objective function in equation (\ref{infinite})
can be combined to accelerate multi-dimensional QTL searches, and especially
permutation testing for significance of an existing candidate QTL.

As indicated in the previous section, equation (\ref{frssx}) can be used
to compute a quantile for the distribution of the residual sum of squares
for a finite size population, i.e. $F_{RSS_x}$, at a specific distance from a hypothetical
QTL explaining all genetic variance. We define a bound for the LogVar objective function based on the quantile of $F_{RSS_x}$, instead of a Lipschitz bound. All boxes are compared against the LogVar transform for a $F_{RSS_x}$ distribution based on the currently found optimum. Pruning of a box is possible if the $1- \epsilon$ quantile for the LogVar distribution at distance $x$ lies below the value evaluated in the centroid of a box with Manhattan radius $x$. If that condition is fulfilled,
then the probability that there would be a new optimum, surpassing the current one, within the box is less than $\epsilon$.

The value for $\epsilon$ will then need to be chosen so that the aggregate
probability of missing a minimum over all splits in a QTL search (and any
permutation testing) is limited to an acceptable level.

In practice, the set of possible box radii appearing in the DIRECT search for a
specific dataset is limited by the structure of the marker map. For a specific minimum
value, the set of quantile limits can therefore be reused after first calculated,
reducing the load of computing the sum of products Guassian-binomial products. Furthermore,
the set of all binomial distribution coefficients for a specific distance can be computed
in almost linear time with respect to the total number of individuals, see Appendix \ref{app2}. If instead each coefficient
would be determined individually when computing the sum, the complexity per term would be
on the order of $O(\log n)$, with a rather larger constant term, even when using an efficient
implementation such as the ones in \citep{boostLib} and \citep{R}.

Finally, the binomial distribution is rather thin-tailed. Therefore, it is not necessary
to sum over all $m_{01}, m_{10}$, but rather only the subset from e.g.\ $np/2 - 8\sqrt{np(1-p)/2}$ to
$np/2 + 8\sqrt{np(1-p)/2}$, still capturing almost all of the probability mass. This reduces the growth
rate of the number of terms in the sum from $O(n^2)$ to $O(n)$. The resulting bound is also only
slightly more conservative, as the quantiles will get shifted upwards. Details for calculating the quantiles can be found in Appendix \ref{app1}.

In all, these adaptations and implementation aspects make it possible to use the defined bound for finite-population online within the Prune-DIRECT algorithm, while avoiding it representing
a major part of the total computational load. Using LogVar as the objective function rather than the
$RSS$ directly is also beneficial even when the bound is not as simple as the Lipschitz bound, since
the expectation will still be a matter of (close to) linear behavior with a bounded slope. This fact
makes the splitting strategy used within DIRECT more efficient.

For intercrosses and other configurations with more than two states, we are still using the finite size
distribution determined for the backcross as an approximation, with an added condition. When more parameters are added,
the computed RSS can fluctuate more. If the RSS and then LogVar value computed in a box is higher than what would
be determined by average in any single linear regression with the same number of degrees of freedom, i.e. the mean of
the $\chi^2$ distribution, then the LogVar transform of the $\chi^2$ mean, rather than the actually computed LogVar value, is
used in determining the pruning condition. At this level small disturbances due to the specific phenotype values overtake
the general Gaussian assumptions.

\section{Results}
The use of DIRECT in \citep{kl3} only resulted in a very small bias when used for permutation testing. Evaluating our PruneDIRECT method with too few permutations and with a specific experimental dataset might render a false positive, as the case where an exhaustive search finds a different optimum can be assumed to be very rare. Also, using a purely simulated dataset might hide non-ideal properties of experimental datasets, giving a validation of our bound which would not work out in practice. For example, the patterns of missing genotype data can give additional confounding and sampling effects.

For these reasons, we decided to use an unpublished experimental dataset with a combination of micro-satellite and SNP markers and varying patterns of missing data for our experiments. The dataset is described in \citep{Nettelblad2009}. Based on inferred haplotypes in the $F_0$ generation, derived using the tool presented in that paper, multiple replicates based on this population structure were constructed and QTL with varying dimensionality simulated accordingly. These replicates could then be analyzed for main effects and permuted runs. By creating hundreds or thousands of replicates, with hundreds of permutations within each, we can expect to discover any deviations between exhaustive search and running PruneDIRECT with a combined termination condition of minimum resolution equivalent to the exhaustive search lattice, and a pruning of impossible split candidates based on the strategy outlined in Section \ref{appl}.

%

\subsection{Simulations}
The simulations were performed on the Tintin cluster at the UPPMAX computational resource center, running as single threads on nodes with AMD Opteron 6220 CPUs. The code is parallel, but since a very large number of replicates were used, many jobs corresponding to individual replicates were executed simultaneously, where each job executed the serial version of the code.

For each run, first the non-permuted QTL model of specified size was fit, allowing all levels of interaction (each genotype-phenotype mean was a free parameter). Then, permutations were created. For exhaustive search, all possible candidate loci sets were explored in a $1$ cM lattice. For PruneDIRECT, the minimum resolution was that same lattice, but in addition the bound was used to avoid splitting of some boxes, if it was definite that the values within that box could not improve on the minimum found in the original main run. The $\epsilon$ value used in PruneDIRECT was $10^{-9}$.

Table \ref{sizes} presents the specific number of replicates, the simulated broad-sense heritability $h^2$, and the number of permutations done for each replicate. Table \ref{times} presents average, minimum and median number of function evaluations for complete sets of main run plus permutations. Timings show that objective function evaluations exceed $90\%$ of the time used in the DIRECT version. The maximum number of function evaluations can exceed the number for an exhaustive search, as the DIRECT implementation we started out from did not natively implement a discretized search grid and thus points coinciding in the discrete lattice could be evaluated multiple times.

 Out of 500 simulated 2D QTL, the loci recovered in 25 of them resulted in more than 10 permuted datasets surpassing the simulated QTL, i.e. less than $99\%$ significant. These runs were also more time-consuming, since the efficiency from the PruneDIRECT algorithm arises from the optimum being rare. If the QTL is not significant, then more boxes will be similar by random and thus fewer boxes can be pruned. If these are removed from the computed number of function evaluations, the average number of evaluations decreases to $7\,529$ thousand. Acceleration compared to exhaustive search is only possible when the result from the main run rises above the noise floor defined by the quantiles and $\chi^2$ distributions. The possible accelerations increase rapidly with stronger QTL signals. Finally, we tested simulating a single four-locus network with a total heritability of $0.30$. Finding this took $9.0$ million function evaluations, which compares favorably to the $963$ millions just finding that network, i.e. omitting any permutation testing, would require using an exhaustive search.

We also created $10,000$ random QTL in a backcross templated on this intercross dataset (by fixing allele origin for either parent). Full concordance was achieved in classical one-dimensional QTL searches between the minimum found through exhaustive search and using PruneDIRECT.


\begin{table}
\caption{\footnotesize\emph{Size of validation tests for two and three dimensions. Time use for validating exhaustive search runs limits the possible size for 3D. Permutations were done per replicate.}}
\vspace{2mm}
\centerline{
\setlength{\leftskip}{0pt}
\label{sizes}
\begin{tabular}{lrrr}
d & Heritability $(h^2)$ & Number of replicates & Number of permutations \\
\hline
2 & 0.09 & 500 & 1000 \\
3 & 0.14 & 500 & 100 \\
\hline
\end{tabular}}
\end{table}
\begin{table}
\caption{\footnotesize\emph{Number of function evaluations for two and three dimensions with exhaustive search and DIRECT, respectively. Minimum, mean, median and maximum number of function evaluations per full run of main QTL search followed by the number of permutations given in Table \ref{sizes} are reported. A very low number of DIRECT runs resulted in full exhaustive searches. Note that the 3-dimensional scan used a lower number of permutations. All numbers are in thousands of function evaluations.}}
\vspace{2mm}
\centerline{
\setlength{\leftskip}{0pt}
\label{times}
\begin{tabular}{rlrrrr}
d & Method & Min & Median & Avg & Maximum \\
\hline
2 \\
& Exhaustive & $76\,126$ & $76\,126$ & $76\,126$ & $76\,126$ \\
& DIRECT & $453$ & $2\,480$ & $11\,012$ & $201\,514$ \\
3 \\
& Exhaustive & $998\,536$ & $998\,536$ & $998\,536$ & $998\,536$ \\
& DIRECT & $1\,565$ & $17\,826$ & $54\,457$ & $323\,628$ \\
\hline
\end{tabular}}
\end{table}
%
\section{Discussion}
This paper presents an efficient scheme, named PruneDIRECT, for simultaneous mapping of several QTL, including permutation testing to determine significance.

The global optimization scheme DIRECT has earlier been adopted to QTL analysis in \citep{kl3}, and the reported speedups are several orders of magnitude \citep{kl3}. The work presented in this paper does not improve on those speedup results. Instead, PruneDIRECT provides an option for performing QTL scans in settings where accuracy and guaranteed results are of of imperative importance. We do so by letting the DIRECT process continue executing until a full exhaustive search has in some sense been performed, but with a pruning taking place removing boxes that are not possible optima. This allows PruneDIRECT to be used not only for finding QTL candidates, but also to be reliably used in permutation testing to find the extreme end of the null hypothesis distribution. These computations can be used to compute significance thresholds as well as assessing the extreme value distribution, something which could form a basis for comparisons between models with different dimensionality and parametrization.

%

We finally propose to use a bound based on quantiles of an analytically derived distribution of the objective function at any distance from a minimum, in order to handle the fact that any real finite-sized population will be affected by the random and discrete nature of recombination, and thus not be perfectly linear. We then used this bound in pruning of the search trees.
 Our new objective function accelerates the performance of the DIRECT process, even when the pruning step is not added. The reason is that DIRECT performs best if the function is linear (finding the top of a single triangle in only a few iterations), and our transformed objective function will in general be almost piecewise linear. This understanding of the expected local form of the objective function, when full genotype information is available, could also be used to better assess probable QTL locations in cases of partial and incomplete genetic information.

It should be noted that the performance for PruneDIRECT is dependent on the heritability. For a trait with no heritability at all, finding the true optimum can in principle only be done by an exhaustive search, as very limited correlations are expected in the objective function between loci. Previous incarnations of DIRECT would have failed in those cases, while PruneDIRECT is adaptive and will perform more function evaluations. If one knows beforehand that QTL with a heritability below some limit $h^2_{min}$ are not relevant, then such information can be added to the pruning process and give better performance even in those cases. We propose that in the future such a scheme could be used to effectively and effortlessly scan for highly significant multi-dimensional QTL in expression QTL data, where tens of thousands of putative phenotypes should be tested.

If the goal is to establish the exact significance level of a highly significant set of QTL, then our PruneDIRECT approach will excel. If most null hypothesis permutations for a dataset have a minimum that is inferior to that determined for the main model, the permuted runs can exit after only a hundred or so function evaluations, even in multiple dimensions. This allows doing runs equivalent to $100,000$ or more permutations for loci where significance levels above $99.9\%$ would be relevant.

It should be noted that our approach for providing accuracy guarantees is based on some underlying assumptions. If phenotype distributions are far from normal, our approximation for the finite-size objective function distribution will not be accurate. It is also necessary that the number of recombinants in any interval is reflected by the mapping distances given. Our approximation takes the randomness of recombination into account, but if a marker map e.g.\ states a zero or very small distance in a region where the actual number of recombinants is much higher, then the approximation breaks down. Hence, it is important to use a genetic map which reflects the actual population studied. Naturally, additional tolerances could be entered into the method by e.g.\ using a higher $p$ value than that provided from the map.  The results presented here could easily be extended to more generations than a $F_2$ intercross or a backcross, as the main difference will be that the rate of crossovers with regard to founder origin is multiplied.


 For linked QTL, PruneDIRECT does well as long as the effects do not cancel
 each other out. When a single QTL effect is fitted the effects from
different QTL are, at least partially, confounded, as there is only a
single variable (the indicator position within the linkage group),
relating back to both components. Among other things, this means that
the total explainable variance can be $0$ at some point in the region
between two linked QTL if the effects at the QTL have opposite
signs. However, outside of the interval between the linked QTL, the behavior is
completely identical to the presence of a single QTL at the position
of the closest QTL in the set, with an effect equivalent to the
combined average effect of all the linked QTL observed from that
position. This can intuitively be understood from the memory-less
nature of the exponential function. To avoid these issues, one could for example imagine enforcing
a coarse-grain splitting of all boxes down to some resolution level and only then applying the prune mode
in PruneDIRECT.

\section*{Acknowledgment}
Per Jensen, Leif Andersson, and Olle K\"ampe are acknowledged for sharing experimental data for evaluation of the method. The computations were performed on resources provided by SNIC through Uppsala Multidisciplinary Center for Advanced Computational Science (UPPMAX).

\section*{Author Disclosure Statement}
No competing financial interests exist.
\newpage
\bibliography{References}
\newpage
\appendix
\appendixpage
\section{Quantile Calculations}\label{app1}
We want to find the value of $rss_x$ such that:
 \begin{align}
 F_{RSS_x}(rss_x|y,Z(0),p)=1-\epsilon
  \end{align}
  One knows that $B(x,n,p)$ for $x$ outside the interval $[\mu-8\sigma,\mu+8\sigma]$ is almost $0$, where $\mu$ and $\sigma$ are the mean and standard deviations of the binomial distribution considered. By using this, one gets:
\begin{align}
&F_{RSS_x}(rss_x|y,Z(0),p)= \\
\nonumber &\sum_{m_{01}'(x)} [B(m_{01}'(x),n_0(0),p)\sum_{m_{10}'(x)} B(m_{10}'(x),n_1(0),p)F_{RSS_x}(rss_x|m_{01}'(x),m_{10}'(x),y,Z(0),p)]+\delta
\end{align}
\begin{align}
\text{    where: } &\delta \approx 0,\\
&  m_{01}'(x),m_{10}'(x)\in \left [\frac{np}{2}-8\sqrt{\frac{np(1-p)}{2}},\frac{np}{2}+8\sqrt{\frac{np(1-p)}{2}} \right ]
\end{align}
An alternate region for  $ m_{01}'(x),m_{10}'(x)$ is to choose $m_{10}'(x)\in[np/2-8\sqrt{np(1-p)/2},np/2+8\sqrt{np(1-p)/2}]$ and then $m_{01}'(x)$ in a circle around $m_{01}'(x)$ with radius $8\sqrt{np(1-p)/2}$. After all, one should find $rss_x$ such that:
\begin{align}
\sum_{m_{01}'(x)} [B(m_{01}'(x),n_0(0),p)\sum_{m_{10}'(x)} B(m_{10}'(x),n_1(0),p)F_{RSS_x}(rss_x|m_{01}'(x),m_{10}'(x),y,Z(0),p)] = 1-\epsilon
\end{align}

This equation can be solved efficiently numerically by e.g. the bisection method. In our implementation, we only solve it down to a resolution level of $0.04$, overestimating the location of the LogVar quantile by at most the equivalent of $1$ cM in the inifite-size population case. Furthermore, since we know the quantile location to be monotonous in terms of increasing $x$, our pre-calculation scheme for different box radii can select the low and high initial bisection bounds based on previously calculated solutions for lower and higher $x$, respectively. In all, this can reduce the amortized number of bisection steps for computing the sum to a handful per unique $x$. For different $x$, the $\phi$ values going into the sum are also unchanged and can thus be stored precalculated, since evaluating the normal CDF is relatively time-consuming.
\section{Calculating the Full Set of Binomial Coefficients} \label{app2}
Most library implementations supporting the binomial distribution compute single values individually. Even if one is using e.g. R \citep{R} and provides a vector as $x$ when computing $B(x, n, p)$, each step will be performed independently. The computational complexity for these steps is non-trivial and on the order of $O(log(x) + log(n))$.

If one knows that all binomial coefficients will be used in a sum, it is instead far more efficient to compute them all as part of the same process. The total computational complexity then becomes linear, with only a few multiplications and divisions per iteration, augmented by a few exponentializations to renormalize the carry values, in order to avoid falling outside the dynamic range of conventional floating point implementations.

\end{document}